\documentclass[a4paper, 12pt]{article}\usepackage[english]{babel}
\usepackage[utf8x]{inputenc}
\usepackage{hyperref, amsmath, amsthm, pgf, latexsym, amsfonts, graphicx, enumerate, float, color, soul, colortbl}
\usepackage{amssymb, url}
\usepackage{hyperref}
\usepackage[affil-it]{authblk}
\usepackage[margin=1in]{geometry}

\usepackage{diagbox}
\usepackage{setspace}

\newtheorem{thm}{Theorem}[section]

\newtheorem{lem}[thm]{Lemma}
\newtheorem{cor}[thm]{Corollary}

\newtheorem{conj}[thm]{Conjecture}

\usepackage{enumitem}
\usepackage{phonetic}

\newcommand{\Z}{\mathbb{Z}}

\newcommand{\R}{\mathbb{R}}

\newcommand{\minz}{\mathrm{min}_0}
\newcommand{\A}{\mathcal{A}}

\newcommand{\B}{\mathcal{B}}

\ExplSyntaxOn
\NewDocumentCommand{\definealphabet}{mmmm}{
    \int_step_inline:nnn{`#3}{`#4}{
        \cs_new_protected:cpx{#1\char_generate:nn{##1}{11}}{\exp_not:N #2{\char_generate:nn{##1}{11}}}
    }
}
\ExplSyntaxOff
\definealphabet{cal}{\mathcal}{A}{Z}
\definealphabet{bb}{\mathbb}{A}{Z}
\definealphabet{scr}{\mathscr}{A}{Z}
\let\>\rangle
\providecommand{\keywords}[1]{\textbf{{ Key Words---}} #1}
\hypersetup{
    colorlinks=true,
    linkcolor=blue,
    filecolor=magenta,      
    urlcolor=cyan,
    pdfpagemode=FullScreen,
    }
\title{Restricted sums of sets of cardinality $2p+1$ in $\Z_p^2$}
\author{Jacinda Eva Terkel\thanks{Email: \href{mailto:jeterkel@memphis.edu}{jeterkel@memphis.edu}}}

\date{}

\begin{document}

\maketitle

% \section{Preface}
% This is a preprint for a paper of mine in additive combinatorics. This copy is for submission to mathematics PhD programs to provide evidence of my research in the field. This submission aims to provide documentation of my interest and experience in researching additive combinatorics.
% If further documentation or questioning are desired regarding my research experience (including a copy of my extensive research notes, or any other relevant documentation) I welcome all requests via e-mail. 

\begin{abstract}
    Let $A\subseteq \Z_p^2$ be a set of size $2p+1$ for prime $p\geq 5$. In this paper, we prove that $A\hat{+}A=\{a_1+a_2\mid a_1,a_2\in A, a_1\neq a_2\}$ has cardinality at least $4p$.  This result is the first advancement in over two decades on a variant of the Erdős-Heilbronn problem studied by Eliahou and Kervaire.
\end{abstract}

\section{Introduction}

In an abelian group $G$ with $A,B\subseteq G$ we write \[A+B=\{a+b\mid a\in A, b\in B\}\] to be the sumset of $A$ and $B$. Similarly, we define  \[A\hat{+}B=\{a+b\mid a\in A, b\in B,a\neq b\}\] to be the \emph{restricted} sumset of $A$ and $B$. Often, we write $2A=A+A$ and $2\hat{\;}A=A\hat{+}A$. A topic of great importance in additive combinatorics is determining the minimum size of $2A$ or $2\hat{\;}A$ given that $A\subseteq G$
 has size $m$. This question has been answered for all abelian groups in the unrestricted case for over twenty years now \cite{SmallSums} but the restricted case remains unsolved in general. More specifically, we are interested in determining the value of the function \[\rho(G,m)=\min\{|2\hat{\;}A|\mid A\subseteq G, |A|=m\}.\] The cases for which the value of this function are known is very limited (see Chapter D.3.1 in \cite{Baj:2018a} for more detail). The case that we are interested in today is a subset of the case where $G$ is an elementary abelian $p$-group for an odd prime $p\geq 5$. This case was first studied extensively by Eliahou and Kervaire and they obtained the following results 
 \begin{thm}[Eliahou and Kervaire  \cite{EKfirst}, \cite{EliKer}]
     If $p\geq 5$ is prime and $p\nmid m-1$ for $m\leq (p^r+3)/2$ \[\rho(\Z_p^r,m)=2m-3.\] 
If $p\mid m-1$ and $m\leq (p^r+3)/2$ we have that
\[2m-3\leq \rho(\Z_p^r,m)\leq 2m-2.\]
If $m=p+1$ we have that
     \[\rho(\Z_p^r,m)=2m-2=2p.\]
 \end{thm}
 This leaves $m=2p+1$ as the smallest unsolved case.
 
In this paper, we determine that $\rho(\Z_{p}^2,2p+1)=4p$ for all prime $p\geq 5$. To do this, we will make use of two of the most ubiquitous theorems in additive combinatorics. Before this, we introduce the notation $\mathrm{min}_0$ which we define as \[\mathrm{min}_0(X)=\max\{0,\min(X)\}.\] ie. $\mathrm{min}_0(S)$ is equal to the minimum of $S$ when $S\subseteq \R_{\geq 0}$ and $0$ if $S$ contains a negative number.

\begin{thm}[Cauchy-Davenport Theorem]
For $A,B\subseteq \Z_p$ for some prime $p$ then we have that \[|A+B|\geq \mathrm{min}\{|A|+|B|-1,p\}.\] \label{cdt}
\end{thm}
% \begin{thm}[Vosper's Theorem]\label{vosper}
% For all non-empty, non-singleton $A,B\subseteq \Z_p$ for prime $p$ with \[p-2\geq |A+B|=|A|+|B|-1\] then $A$ and $B$ are arithmetic progressions with common difference.
% \end{thm}
\begin{thm}[Dias Da Silva and Hamidoune \cite{DDSandH}]
\label{ehc}
If $p$ is prime with $A,B\subseteq \Z_p$ we have that \[|A\hat{+}B|\geq\mathrm{min}_0\{|A|+|B|-3,p\}.\] Furthermore, if $|A|\neq |B|$ then
\[|A\hat{+}B|\geq\mathrm{min}_0\{|A|+|B|-2,p\}.\] 
\end{thm}
Our main result is as follows:
\begin{thm}\label{maintheorem}
If $p\geq 5$ is prime then \[\rho(\Z_p^2,2p+1)=4p.\]
\end{thm}

We first establish this simple yet consequential result.

\begin{lem}\label{translation}
For $g\in \Z_p^2$ and $A\subseteq \Z_p^2$ we have that
    \[|2\hat{\;}A|=|2\hat{\;}(A+g)|.\]
\end{lem}
\begin{proof}
    Note that $a_1=a_2\iff a_1+g=a_2+g$, and so we have that $2\hat{\;}(A+g)=2g+2\hat{\;}A$, from which our claim follows.
\end{proof}
Let $H$ be some non-trivial proper subgroup of $\Z_p^2$ (ie. $H\cong\Z_p$). Similarly, index the cosets of $H$ by $H_0,H_1,\dots,H_{p-1}$ where $H_0=H$ and $H_{i}+H_{j}=H_{i+j}$, and index their intersections with $A$ such that $A_i=A\cap H_i$ and similarly let $B_i=2\hat{\;}A\cap H_i$. Sometimes we may refer to these indexed sets by indices outside of the range $[0,p-1]$, and these should be identified with their representative within the aforementioned interval modulo $p$ (ie. $H_i=H_{kp+i}$ for all $k\in\Z$ and $i\in [0,p-1]$).

It is obvious that \[|2\hat{\;}A|=\sum_{i=0}^{p-1}|B_i|\] as $B_0,B_1,\dots,B_{p-1}$ is a partition of $2\hat{\;}A$. Additionally, we see that \[B_i=\bigcup_{j=0}^{p-1}(A_{j}\hat{+}A_{i-j}),\] and so \[|B_i|\geq\max_{0\leq j\leq p-1}|A_{j}\hat{+}A_{i-j}|.\] 
Note that when $j\neq i-j$, $A_j$ and $A_{i-j}$ are disjoint and so $A_{j}\hat{+}A_{i-j}=A_{j}+A_{i-j}$. Thus, by Theorem \ref{cdt}, Theorem \ref{ehc}, and Lemma \ref{translation}  we have that \begin{equation}\label{CosetLBMain}
    |B_i|\geq \max_{0\leq j\leq p-1}\{\mathrm{min}_0\{|A_j|+|A_{i-j}|-1-2\epsilon_{i,j},p\}\mid \emptyset\not\in \{A_j,A_{i-j} \}\}
\end{equation}
where $\epsilon_{i,j}=1$ if $j=i-j$ and $\epsilon_{i,j}=0$ otherwise.
Now, let $\alpha$ be some index for which $|A_\alpha|\geq |A_i|$ for all $i$.  By the pigeonhole principle, it is seen that $|A_\alpha|\geq 3$. By Lemma \ref{translation}, we may assume, without sacrificing generality, that $\alpha=0$ Additionally, let $m$ be the number of non-zero $i$ for which $A_i$ is non-empty and let $S=\{1\leq i\leq p\mid A_i\neq \emptyset\}.$  Similarly, let $T=\{1\leq i\leq p\mid B_i\neq \emptyset\}$, $S_0=S\cup\{0\}$ and $T_0=T\cup\{0\}$.

Observe that 
\begin{equation}\label{basiclb}
    |2\hat{\;}A|=\sum_{i=0}^{p-1}|B_i|\geq |2\hat{\;}A_0|+\sum_{s\in S}|A_0+A_s|.
\end{equation}
and since $|A_s|\leq |A_0|$ if $|A_0|\leq (p+1)/2$ then we have that \[ |2\hat{\;}A|\geq |2\hat{\;}A_0|+\sum_{s\in S}(|A_0|+|A_s|-1)\] which can be taken advantage of in a multitude of ways. For this reason, we consider the cases of $|A_0|\leq (p+1)/2$ and $|A_0|\geq (p+3)/2$ separately Specifically, we define \textbf{Case 1:} $|A_0|\leq \frac{p+1}{2}$ and  \textbf{Case 2:} $|A_0|\geq \frac{p+3}{2}$.

\section{Case 1}

As mentioned before, the reason we divided this problem into cases based on $|A_0|$ is that, no matter what, in Case 1 we have that $|A_i|+|A_j|-1\leq p$, and so \[|A_i+A_j|\geq \minz\{p,|A_i|+|A_j|-1-2\epsilon_{i,j}\}=|A_i|+|A_j|-1-2\epsilon_{i,j}\] as per (\ref{CosetLBMain}).

We will make use of the following definitions:
\[\A=\{A_i\mid |A_i|\neq \emptyset\},\]
\[\A'=\{A_i\mid |A_i|\neq \emptyset,i\neq0\}\]
\[\A_w=\{A_i\mid |A_i|=w\},\]
\[\B_w=\{B_i\mid |B_i|=w\},\]
% \[\A_w'=\{A_i\mid |A_i|\geq w\},\]
% \[\B_w'=\{B_i\mid |B_i|\geq w\},\]
\[C_w=|A_w|,\quad\quad D_w=|\B_w|,\quad\quad D_w'=\sum_{i\geq w}|\B_i|.\]
It is clear that \[|A|=\sum_{w=1}^{p}wC_w\quad\quad\text{and}\quad\quad|2\hat{\;}A|=\sum_{w=1}^{p}wD_w.\] Also observe that 
Now, note that $D_{w}'-D_{w+1}'=D_w$, and so we have that \begin{equation}
    \label{dSum}|2\hat{\;}A|=\sum_{w=1}^{p}wD_w=\sum_{w=1}^{p}w(D_{w}'-D_{w+1}')=\sum_{w=1}^{p}D_w'.
\end{equation}

Let $m=|S|=|\mathcal{A}|-1$ (ie. $m$ is the number of non-zero $i$ for which $A_i$ is non-empty.)

By Theorem \ref{ehc}, there exists at least $\minz\{2|\A|-3,p\}=\minz\{2m-1,p\}$ distinct values of $i$ for which there is some $X\in 2\hat{\;}\A$\footnote{Because the elements of $\A$ are pairwise disjoint and non-empty: every element of $2\hat{\;}\A$ is non-empty.} satisfying $X\subseteq H_i$. In other words, $D_1'\geq \minz\{2m-1,p\}$. This is the basis for our first instance of subcases, those being

\begin{enumerate}[label=(1\Alph*)]
    \item $m\geq \frac{p+1}{2}$ (in which case $B_i$ must be non-empty for all $i$);
    \item $m\leq \frac{p-1}{2}$ (in which case we must have that $|T\setminus S|\geq m-2$).
\end{enumerate}

Before entering the subcases, we prove the following statement regarding Case 1 in general
\begin{lem}\label{Case1Backbone}
In Case 1
\[|2\hat{\;}A|\geq (m+1)(|A_0|-1)+2p-1+\sum_{i\in T\setminus S}|B_i|.\]
\end{lem}
\begin{proof}
Using (\ref{CosetLBMain}), Theorem \ref{cdt}, and Theorem \ref{ehc} and the fact that $|A_i|+|A_j|-1\leq 2|A_0|-1\leq p$ for all $i,j$\footnote{From the fact that we are in Case 1.} we get that
    \[\sum_{i\in S_0}|B_i|\geq 2|A_0|-3+\sum_{i\in S}\left(|A_0|+|A_i|-1\right)=(m+2)|A_0|-3-m+\sum_{i\in S}|A_i|\]\[=(m+2)(|A_0|-1)-1+\sum_{i\in S}|A_i|=(m+2)(|A_0|-1)-1+(|A|- |A_0|)\]\[=(m+1)(|A_0|-1)+2p-1.\]

Since\[|2\hat{\;}A|=\sum_{B_i\in \B}|B_i|\] we have that 
\[|2\hat{\;}A|=\sum_{i\in T_0}|B_i|=\sum_{i\in S_0}|B_i|+\sum_{i\in T\setminus S}|B_i|\geq (m+1)(|A_0|-1)+2p-1+\sum_{i\in T\setminus S}|B_i|,\] and so our claim is proven.
\end{proof}

\subsection{Case 1A}

In Case 1A, we add the additional assumption that $m\geq \frac{p+1}{2}$.

\begin{lem}
    In Case 1A, if $|2\hat{\;}A|<4p$ then $|A_0|= 3$.
\end{lem}
\begin{proof}
    As mentioned in the introduction, $|A_0|\geq 3$ by the pigeonhole principle, and so it suffices to prove that $|A_0|\leq 3$. 
    Indeed, since $|B_i|\geq 1$ for all $i$ (implying $|T|=p-1$) with the help of Lemma \ref{Case1Backbone}  we get that \[|2\hat{\;}A|\geq (m+1)(|A_0|-1)+2p-1+\sum_{i\in T\setminus S}|B_i|\geq (m+1)(|A_0|-1)+2p-1+(p-|\A|)\] 
    \[=3p+(m+1)(|A_0|-2)-1.\] If indeed it is true that $|2\hat{\;}A|<4p$ then the above implies \[4p-1\geq 3p+(m+1)(|A_0|-2)-1\] which then gives us that \[ |A_0|\leq 2+\frac{p}{m+1}\leq 2+\frac{2p}{p+3}=4-\frac{2}{p+3}<4.\] Thus, $|A_0|\leq 3$ as it is an integer.
\end{proof}

From this, it follows that every set in Case 1A satisfies \begin{equation}\label{Csum1}
    C_0+C_1+C_2+C_3=p,
\end{equation}
\begin{equation}\label{Csum2}
    3C_3+2C_2+C_1=2p+1,
\end{equation} and \begin{equation}\label{Csum3}
    C_1+C_2+C_3=m+1.
\end{equation}

% It is already known that $C_3\geq 1$, and we show a similar yet weaker (albeit still useful) fact for $C_1$ and $C_2$
% \begin{lem}\label{Case1ANoHomo}
% In Case 1A if $|2\hat{\;}A|\leq 4p$ then
%     \[C_1+C_2\geq 1.\]
% \end{lem}
% \begin{proof}
%     Note that if the claim does not hold, by (\ref{Csum2}) we must have that \[C_3=\frac{2p+1}{3}.\]  Assume this is true. Since we are in Case 1A, it follows that, for all $i$, there is some $X\in2\hat{\;}\A$ where $X\subseteq H_i$. Additionally, since $X=A_{i_1}+A_{i_2}$ for $i_1\neq i_2$ and $A_{i_1},A_{i_2}$ non-empty. It follows from Theorem \ref{cdt} that $|X|\geq |A_{i_1}|+|A_{i_2}|-1=5$, and thus $|B_j|\geq 5$ for all $j$ implying $|2\hat{\;}A|\geq 5p>4p$.
% \end{proof}

\begin{lem}\label{Case1A:D3LB}
  In Case 1A \[D_3'\geq \minz\{2C_3+C_2+C_1-1,p\}.\] 
\end{lem}
\begin{proof}
    Via Theorem \ref{cdt}, we have that \[\A+\A_3=\{A_i\hat{+}A_j\mid A_i\in \A, A_j\in \A_3\}\] contains members which are subsets of at least $\minz\{|\A_3|+|\A|-1 ,p\}=\minz\{2C_3+C_2+C_1-1,p\}$ distinct cosets. Additionally, via (\ref{CosetLBMain}) we have that $|X|\geq 3$.  for all $X\in \A+\A_3$, and from this our claim follows\footnote{While it is true that the elements $A_i\hat{+}A_j$ of $\A+\A_3$ are restricted sums, for $(A_i,A_j)\in \A\times \A_3$ we have that $A_i\hat{+}A_j=A_i+A_j$ when $|A_i|\neq 3$ implying $|A_i\hat{+}A_j|\geq |A_i|+|A_j|-1\geq |A_j|=3$, and if $|A_i|=3$ then $|A_i\hat{+}A_j|\geq |A_i|+|A_j|-3=3$.}.
\end{proof}
Similar arguments on \[\A_3+\A_2=\{A_i\hat{+}A_j\mid A_i\in \A_2, A_j\in \A_3\}=\{A_i+A_j\mid A_i\in \A_2, A_j\in \A_3\}\] and \[2\hat{\;}\A_3=\{A_i\hat{+}A_j\mid A_i,A_j\in \A_3, i\neq j\}=\{A_i+A_j\mid A_i,A_j\in \A_3, i\neq j\}\] result in
\begin{lem}\label{Case1A:D4LB}
    In Case 1A
    \[D_4'\geq \minz\{C_3+C_2-1,p\}.\]
\end{lem}
and
\begin{lem}\label{Case1A:D5LB}
    In Case 1A
    \[D_5'\geq \minz\{2C_3-3,p\}.\]
\end{lem}
respectively.

  Now, we utilize these Lemmas.
\begin{lem}\label{Case1A:D3lp}
        In Case 1A, if $|2\hat{\;}A|<4p$ and $D_3\neq p$ then $3\geq 2C_3+C_1+C_0$
\end{lem}
\begin{proof}
   Since $D_3'\neq p$, it follows that $D_4',D_5'\neq p$ either. Together with facts that $D_1'=p$, $D_2'\geq D_3'$,  Lemma \ref{Case1A:D3LB}, Lemma \ref{Case1A:D4LB}, Lemma \ref{Case1A:D5LB}, and (\ref{dSum}) we get that
\[|2\hat{\;}A|\geq D_5'+D_4'+D_3'+D_2'+D_1'\geq p+D_5'+D_4'+2D_3'\geq p+ 7C_3+3C_2+2C_1-6.\] With (\ref{Csum2}) and our assumption $|2\hat{\;}A|\leq 4p-1$ we get that \[4p-1\geq 5p-4+C_3-C_2,\] and via some rearrangement we get that \[C_2+3\geq p+C_3,\] and via substitution of (\ref{Csum1}) this gives us\[3\geq 2C_3+C_1+C_0\]
\end{proof}
% \begin{lem}\label{Case1A:D3p}
%       In Case 1A, if $|2\hat{\;}A|<4p$ and  $D_3'=p$ then $p-3\leq m$.
% \end{lem}
% \begin{proof}
% If this is the case, then by Lemma \ref{Case1Backbone} we have that \[|2\hat{\;}A|\geq (m+1)(|A_0|-1)+2p-3+3(p-m-1)=(m+1)(|A_0|-4)+5p-3.\] 
% \[=5p-3-(m+1).\] Thus, we must have that \[5p-3-(m+1)\leq 4p-1\] implying \[p-3\leq m\]
% \end{proof}
\begin{lem}\label{Case1A:D3p}
    In Case 1A, if $|2\hat{\;}A|<4p$ and $D_3'=p$ then $C_0+C_3\leq 2$.
\end{lem}
\begin{proof}
Clearly, if $D_4'=p$ or $D_5'=p$ then $|2\hat{\;}A|\geq 4p$, and so we may assume $D_4'$ and $D_5'$ are both less than $p$, and so by Lemma \ref{Case1A:D4LB}, Lemma \ref{Case1A:D5LB}, and (\ref{dSum}) we have that 
\[|2\hat{\;}A|\geq D_5+D_4+3p\geq 3C_3+C_2-4+3p.\] If we have that $4p-1\geq  3C_3+C_2-4+3p$ which implies \[p+3\geq 3C_3+C_2.\] By (\ref{Csum2}) we now have that \[C_1+C_2\geq p-2,\] and thus $C_0+C_3\leq 2$ by (\ref{Csum1}).
\end{proof}
\begin{lem}\label{Case1A:C3e1}
    In Case 1A if $C_3=1$ then $C_0=C_1=0$.
\end{lem}
\begin{proof}
   If $C_3=1$ then by (\ref{Csum1}) and (\ref{Csum2}) we have that \[2p-2=2C_2+C_1\leq 2C_2+2C_1+2C_0=2p-2\] implying $C_1+2C_0=0$ and our claim follows.
\end{proof}
\begin{lem}\label{Case1A:c3e2}
    In Case 1A, if $C_3=2$, $C_0=0$ then $|2\hat{\;}A|\geq 4p$.
\end{lem}
\begin{proof}
    In the case of $C_3=2$ and $C_0=0$, by (\ref{Csum2}) and (\ref{Csum1}) we have that $C_2=p-3$ and so $C_1=1$. Without loss of generality, let $|A_0|=3$, and also $|A_x|=3$, and $|A_y|=1$ for $x\neq 0$. For all other $i$ we have that $|A_i|=2$. (\ref{CosetLBMain}) now gives us that \[|B_i|\geq \max_{0\leq j\leq p-1}\{\mathrm{min}_0\{|A_j|+|A_{i-j}|-1-2\epsilon_{i,j},p\}\mid \emptyset\not\in \{A_j,A_{i-j} \}\}\]
    which implies
    \[|B_i|\geq \max\{\minz\{p,|A_x|+|A_{i-x}|-1\},\minz\{p,|A_0|+|A_{i}|-1\}\},\] and because $|A_x|=|A_0|=3$, and $p\geq 5$ we have that \[|B_i|\geq 2+\max\{|A_{i-x}|,|A_{i}|\}.\]  But now, since $x\neq 0$ it follows that $A_{x}\neq A_{i-x}$, and thus it follows that at least one of $A_x$ and $A_{i-x}$ are not $y$ which then means that at least one of them has cardinality greater than or equal to $2$ meaning for all $i$ we have that $|B_i|\geq 4$ implying $D_4'=p$ and thus $|2\hat{\;}A|\geq 4p$.

    By Lemma \ref{Case1A:D5LB} we have that $D_5'\geq 1$.
\end{proof}
Combining Lemma \ref{Case1A:D3lp}, Lemma \ref{Case1A:D3p}, Lemma \ref{Case1A:C3e1}, and Lemma \ref{Case1A:c3e2} we can now do the following.
\begin{cor}
In Case 1A, if $|2\hat{\;}A|<4p$ then $C_3=1$, $C_2=p-1$, $C_1=0$, and $C_0=0$.
\end{cor}
\begin{proof}
Keep in mind throughout this proof that $C_3\geq 1$ by the pigeonhole principle.

If $D_3'=p$ then by Lemma \ref{Case1A:D3p} we have that $C_3+C_0\leq 2$. Thus $C_3\leq 2$. Since $C_3\neq 0$ we have either $C_3=2$, in which case $C_0=0$ by the inequality or $C_3=1$ which by Lemma \ref{Case1A:C3e1} implies that $C_0=C_1=0$.

In the case of $C_3=2$ and $C_0=0$, Lemma \ref{Case1A:c3e2} implies that $|2\hat{\;}A|\geq 4p$.

If $D_3'\neq p$ then Lemma \ref{Case1A:D3lp} $3\geq 2C_3+C_1+C_0 $ directly implies $C_3=1$, and so by Lemma \ref{Case1A:C3e1} we have that $C_0=C_1=0$.

Thus, regardless of $D_3'$ if $|2\hat{\;}A|<4p$ then we must have that $C_3=1$, $C_1=0$, $C_0=0$, and thus via extension by (\ref{Csum1}): $C_2=p-1$.

%     Consider first if $D_3'=p$ in which case we refer to Lemma \ref{Case1A:D3p} to get that $C_0+C_3\leq 2$. However, recall that $C_3\geq 1$ by the pigeonhole principle and so $C_0\in\{0,1\}$. 
    
%     If $C_0=1$ then it is implied that $C_3=1$ and so via (\ref{Csum2}) we get that \begin{equation}\label{Case1A:corEq1}
%         C_1+2C_2=2p-2
%     \end{equation} However, since $C_0+C_3=2$ we have that $C_1+C_2=p-2$ by (\ref{Csum1}) which then implies that $C_2=p$ by (\ref{Case1A:corEq1}), which cannot be, and so, we have that $C_0=0$ implying $m+1=p$.

% Similarly, with Lemma \ref{Case1A:D3p} giving us that \[3\geq 2C_3+C_1+C_0\] we have that $C_3\leq 1$, and since (as mentioned before) $C_3\geq 1$ we have that $C_3=1$ and so $C_1+C_0\leq 1$. If $C_0=0$ then we use the same argument as above to demonstrate that we cannot have $(C_0,C_3)$
    \end{proof}

With the potential number of $2p+1$-subsets narrowed down in Case 1A drastically, we ask the reader to note that if $|2\hat{\;}A|<4p$ then (WLOG via Lemma \ref{translation}) we have that $|A_0|=3$ and $|A_i|=2$ for all non-zero $i$.

Now we must utilize another famous addition theorem.

\begin{thm}[Vosper \cite{Vosper}] \label{vosper}
  If $A,B\subseteq \Z_p$ satisfy $2\leq |A|,|B|$ then \[|A+B|\leq \min\{|A|+|B|-1,p-2\}\] if and only if $A$ and $B$ are arithmetic progressions with a common difference.
\end{thm}

From this, we can prove the following

\begin{lem}\label{Case1A:arith}
    In Case 1A for $p\geq 7$, if $|2\hat{\;}A|<4p$, then there exists some $d$ for which each $A_i$ is an arithmetic progression of difference $d$ and this $d$ is the same for all $A_i$.
\end{lem}
\begin{proof}
    First, note that if $A_0$ was not an arithmetic progression then by Theorem \ref{vosper}, for $p\geq 7$, we would have that \[|B_i|\geq |A_0+A_i|\geq |A_0|+|A_i|\geq 5\] for all non-zero $i$, and because $|B_0|\geq |2\hat{\;}A_0|\geq 2|A_0|-3=3$ we would have that $|2\hat{\;}A|\geq 5p-2> 4p$. 

    Thus, $A_0$ is an arithmetic progression with some difference $d$.

    Since for non-zero $i$ we have $|A_i|=2$ it is trivial that $A_i$ is an arithmetic progression (let us say with difference $d_i$). But I now claim that if $|2\hat{\;}A|\leq 4p-1$ then for all $i$: $d_i=d$. We can prove this as follows: By Theorem \ref{vosper} observe that for $i\neq 0$ we have  $|B_i|\geq |A_0|+|A_i|-1=5-\epsilon_i$ where $\epsilon_i=\begin{cases}
        1 & d_i=d\\ 0 & d_i\neq d
        
    \end{cases}$. Recall that $|B_0|\geq 3$, and by Assuming $|2\hat{\;}A|\leq 4p-1$ and letting $E$ be the number of non-zero $i$ for which $\epsilon_i=1$ we have that
    \[4p-1\geq 5(p-1)-E+3=5p-2-E\geq 4p-1,\] and so we equality holds throughout implying $E=p-1$ and our claim is proven.
\end{proof}

By Lemma \ref{translation} it suffices to consider only when $A_0$ takes the form \[A_0=\{0,d,2d\}\] for some non-zero $d\in H$. Additionally, for non-zero $i$ define $a_i$ such that \[A_i=\{a_i,a_i+d\}.\] It should be noted that by (\ref{CosetLBMain}) that \begin{equation}\label{Case1A:LBplan}
    |B_i|\geq\begin{cases}
    3 & i=0;\\ 4 & i\neq 0.
\end{cases},
\end{equation}  and so if we are to have $|2\hat{\;}A|\leq 4p-1$ then equality must hold in (\ref{Case1A:LBplan}) for all $i$.

We now prove some facts regarding our $a_i$.

\begin{lem}
\label{Case1A:BuildFor}
    In Case 1A for $p\geq 7$, if $|2\hat{\;}A|\leq 4p-1$ then all of the following hold for non-zero $i,j$ with $i\neq j$:
    \begin{enumerate}
        \item $B_i=\{a_i,a_i+d,a_i+2d,a_i+3d\}$,
        \item $a_j+a_{i-j}\in A_i=\{a_i,a_i+d\}$,
        \item $a_{2i}=2a_i+\delta d$ for some $\delta\in\{-2,-1,0,1\}$,
        \item $B_0=\{d,2d,3d\}$,
        \item $a_i+a_{-i}=d$.
    \end{enumerate}
\end{lem}
 \begin{proof}
With Lemma \ref{Case1A:arith} and the above discussion in mind we can prove the statements as follows: 
 
     The first claim follows from the fact that we must have that $|B_i|=4$, $A_0+A_i=\{a_i,a_i+d,a_i+2d,a_i+3d\}$ has size $4$ and $A_0+A_i\subseteq B_i$ meaning $A_0+A_i=B_i$.

     To prove the second claim we see that $A_{j}+A_{j-i}=\{a_{j}+a_{i-j},a_{j}+a_{i-j}+d,a_{j}+a_{i-j}+2d\}\subseteq B_i$ and since $p\geq 7$ this means we must either have that $a_{j}+a_{i-j}=a_i$ or $a_{j}+a_{i-j}=a_i+d$, ie. $a_{j}+a_{i-j}\in A_i$.

     For the third claim, we similarly observe that $2\hat{\;}A_{i}=\{2a_i+d\}\subseteq B_{2i}=\{a_{2i},a_{2i}+d,a_{2i}+2d,a_{2i}+3d\}$ and we see our claim follows.

     For the fourth claim follows from the facts that $|B_0|=3$ and $2\hat{\;}A=\{d,2d,3d\}\subseteq B_0$ like our proof of Claim 1.

     For the fifth claim we see that $A_i+A_{-i}=\{a_i+a_{-i},a_i+a_{-i}+d,a_i+a_{-i}+2d\}\subseteq B_0=\{d,2d,3d\}$ implying that $a_i+a_{-i}=d$. 
 \end{proof}

Let us define $\mu_i=\frac{a_{i}-(a_{i-1}+a_1)}{d}$. By Lemma \ref{Case1A:BuildFor}, we have that $\mu_{i}\in\{0,1\}$ for all $i\in[3,p-1]$ and $\mu_2\in\{-2,-1,0,1\}$.

Thus, we have the recurrence relation $a_{i+1}=a_1+a_i+\mu_i$ based on a predefined $a_i$ which gives us   \[d-a_1=a_{-1}=a_{p-1}=a_1+\sum_{i=2}^{p-1}(a_1+d\mu_i)=(p-1)a_1+d\sum_{i=2}^{p-1}\mu_i,\] and so it is implied that \[\sum_{i=2}^{p-1}\mu_i=1 \mod p,\] and because the sum cannot exceed $p$ or go below $-1$, the implication is that $\sum_{i=2}^{p-1}\mu_i=1$ exactly.

Because $\mu_i\geq 0$ for $i\neq 2$ and $\mu_2\in\{-2,-1,0,1\}$ we have that the number of $i$ (other than $2$) for which $\mu_i=1$ is $1-\mu_2\in\{0,1,2,3\}$.

Thus, it follows that, for every $a_i$ we have that for some $\mu_2 \leq u\leq 1-\mu_2$ we have that $a_i-ud\in K=\<a_1\>$ Thus, by the 1st and 4th statements in Lemma \ref{Case1A:BuildFor} we have that for any $a\in A$ we there exists an integer $u$ within satisfying $\mu_2\leq u\leq 4-\mu_2$ such that $a-ud\in K$. However, this then implies that there are at most $5$ cosets $K_i$ of $K$ where $K_i\cap A$ is non-empty, and so we have that there is some other subgroup of $\Z_p^2$ that intersects $A$ at most $5\leq \frac{p-1}{2}$ different cosets and so we have as follows:
\begin{lem}\label{Case1A:Reduction}
    In Case 1A, for $p\geq 11$, if $|2\hat{\;}A|<4p$ then there is an instance in Case 1B or Case 2 with $|2\hat{\;}A|<4p$. This then implies that if one manages to prove that $|2\hat{\;}A|\geq 4p$ in Case 1B and Case 2 then $|2\hat{\;}A|\geq 4p$ in Case 1A.
\end{lem}
With this, we move towards proving that $|2\hat{\;}A|$ in Case 1B and Case 2.

\subsection{Case 1B}

In Case 1B, we assume that we are not in Case 1A but we are still in Case 1 meaning $m\leq \frac{p-1}{2}$.

\begin{lem}\label{Case1B:dbound}
    In Case 1B, if $|2\hat{\;}A|\leq 4p-1$ then $|A_0|(m+1)-|A|\leq 2$.
\end{lem}
\begin{proof}
    In this case, since $|\B|\geq \minz\{p,2|\A|-3\}=\minz\{p,2m-1\}$, we guarantee the existence of at least $m-2$ distinct $i$ such that $i\in T\setminus S$. Let $d=|A_0|(m+1)-|A|$.  Now, by Lemma \ref{Case1Backbone} we have that \[4p-1\geq |2\hat{\;}A|\geq (m+1)(|A_0|-1)+2p-1+\sum_{i\in T\setminus S}|B_i|\geq (m+1)(|A_0|-1)+2p-1+(m-2)\] implying
    \[2p+3\geq |A|+d=2p+1+d,\] and our claim follows.
\end{proof}

\begin{lem}
    In Case 1B, $|2\hat{\;}A|\geq 4p$.
\end{lem}
\begin{proof}
    Lemma \ref{Case1B:dbound} implies that there exists some selection of $\omega,\psi\in S$ such that for all $i\in S'=S_0\setminus\{\omega,\psi\}$ we have that $|A_i|=|A_0|$, and also $2|A_0|-2\leq |A_\psi|+|A_\omega|\leq 2|A_0|$.

    From this and (\ref{CosetLBMain}), we may deduce that \begin{equation}
    \label{Case1B:ineq}
        |B_i|\geq \begin{cases}
        2|A_\omega|-3 & i=2\omega;\\
        2|A_\psi|-3 & i=2\psi;\\
        2|A_0|-3 & \text{otherwise}.\\
    \end{cases}
    \end{equation}
Additionally, via \ref{cdt} there must exist at least $\minz\{2|\A|-1,p\}=2m+1$ distinct $x\in[0,p-1]$ such that $x=i+j$ for some (not necessarily distinct) $i,j\in S_0$ Let the set of such $x$'s be $\mathcal{X}$. We account for $m+1$ of these via $0+i=i$ for $i\in S_0$, and so using Lemma \ref{Case1Backbone}, (\ref{Case1B:ineq}), the facts that $(m+1)|A_0|=2p+1+d$, $|A_\omega|+|A_\psi|\geq 2|A_0|-d$, $|A_0|\leq (p+1)/2$, and $m\leq (p-1)/2$ we have that \[|2\hat{\;}A|\geq (m+1)(|A_0|-1)+2p-1+\sum_{i\in \mathcal{X}\setminus S_0}|B_i|\]
\[\geq (m+1)(|A_0|-1)+2p-1+(m-2)(2|A_0|-3)+(2|A_\omega|-3)+(2|A_\psi|-3)\]
\[=4p+d-(m+1)+(m-2)(2|A_0|-3)+(2|A_\omega|-3)+(2|A_\psi|-3)\]
\[=4p+d-m-1+(2m|A_0|-4|A_0|-3m+6)+(4|A_0|-2d-6)\]
\[=4p+2m|A_0|-4m-d-1\]\[=4p+2(m+1)|A_0|-4m-d-1-2|A_0|\]
\[=8p+d+1-2|A_0|-4m\geq 8p+1-(p+1)-2(p-1)=5p+2>4p.\]

Our claim now follows from the above and Lemma \ref{Case1B:dbound}.  
\end{proof}

\section{Case 2}
\label{Section2}

In this section, we yet again introduce more terminology. Let $\ell$
be the number of non-zero $i$ for which $|A_0|+|A_i|-1\geq p$, and let $s$ be the number of non-zero $i$ where $A_i$ is non-empty and $|A_0|+|A_i|-1< p$. It follows that $m=\ell+s$. This distinction is made as $\ell$ is the number of $i\in S$ for which $|A_0+A_i|$ is guaranteed to have size $p$ per Theorem \ref{cdt}.
We now will move towards proving a Lemma akin to Lemma \ref{Case1Backbone}, but instead for Case 2.

\begin{lem}\label{Case2:2A0}
    In Case 2, $|2\hat{\;}A_0|=p$.
\end{lem}
\begin{proof}
    From Theorem \ref{ehc}, we have that \[p\geq |2\hat{\;}A_0|\geq \minz\{2|A_0|-3,p\}\geq \minz\left\{2\frac{p+3}{2}-3,p\right\}=\minz\{p,p\}=p.\]
\end{proof}

\begin{lem}\label{Case2Backbone}
    In Case 2, \[|2\hat{\;}A|\geq (l+1)p+s|A_0|+\sum_{i\in T\setminus S}|B_i|.\]
\end{lem}
\begin{proof}
    Keeping (\ref{CosetLBMain}) and specifically Theorem \ref{cdt} in mind we have that
     \[|2\hat{\;}A|= \sum_{i\in T_0}|B_i|=\sum_{i\in S_0}|B_i|+\sum_{i\in T\setminus S}|B_i|\geq |2\hat{\;}A_0|+\sum_{i\in S}|B_i|+\sum_{i\in T\setminus S}|B_i|,\] and now with Lemma \ref{Case2:2A0} we get that \[|2\hat{\;}A|\geq p+\sum_{i\in S}|B_i|+\sum_{i\in T\setminus S}|B_i|\geq p+\sum_{i\in S}|A_0+A_i|+\sum_{i\in T\setminus S}|B_i|\]
     \[\geq (l+1)p+s|A_0|+\sum_{i\in T\setminus S}|B_i|.\]
\end{proof}

We now will demonstrate that, for each value of $\ell$ we have that $|4\hat{\;}A|\geq 4p$

\begin{lem}\label{Case2:lg3}
In Case 2, If $\ell\geq 3$ then $|2\hat{\;}A|\geq 4p$.
\end{lem}
\begin{proof}
    From Lemma \ref{Case2Backbone}, if $\ell\geq 3$ then \[|2\hat{\;}A|\geq 4p+s|A_0|+\sum_{i\in T\setminus S}|B_i|\geq 4p. \]
\end{proof}
For $\ell \leq 2$, we must often provide special consideration to smaller values of $s$.
\begin{lem}\label{Case2:L2S2}
    In Case 2, If $\ell=2$ and $s\geq 2$ then $|2\hat{\;}A|\geq 4p$.
\end{lem}\
\begin{proof}
    From Lemma \ref{Case2Backbone}, if $\ell=2$ then \[|2\hat{\;}A|\geq 3p+s|A_0|+\sum_{i\in T\setminus S}|B_i|\geq 3p+s\frac{p+3}{2},\] and so if $s\geq 2$ then $|2\hat{\;}A|\geq 4p$.
\end{proof}

\begin{lem}\label{Case2:L2S1}
    In Case 2, If $\ell=2$ and $s=1$ then $|4\hat{\;}A|\geq 4p$.
\end{lem}
\begin{proof}
     In this case we may define $\beta,\gamma,$ and $\delta$ to be the three distinct elements of $[1,p-1]$ such that $A_\beta,A_\gamma,$ and $A_\delta$ are not empty satisfying \[|A_0|\geq |A_\beta|\geq |A_\gamma|\geq p+1-|A_0|>|A_\delta|\]
    and \[|A_0|+|A_\beta|+|A_\gamma|+|A_\delta|=2p+1.\]
    These conditions intersect to give us that \[2p+1<p+1+|A_\beta|+|A_\gamma|\] or \[p<|A_\beta|+|A_\gamma|,\] and thus $|A_\beta|\geq \frac{p+1}{2}$.

    Now see that \[S=\{0,\beta,\gamma,\delta\}\] and \[\beta+S=\{\beta,2\beta,\beta+\gamma,\beta+\delta\}\] must not be the same set as this would imply that the sets have that same sum, and thus \[4\beta=0\] which cannot be as $\beta\neq 0$. Thus, for some $\iota\in S$ we have that $\beta+\iota\in T\setminus S$. If $\iota=\beta$ then it is seen that \[|B_{\beta+\iota}|=|B_{2\beta}|\geq|2\hat{\;}A_\beta|\geq2|A_\beta|-3\geq |A_\beta| +\frac{p-5}{2}\geq  |A_\beta|\] as $p\geq 5$. It is also observed that if $\iota\neq \beta$ then \[|B_{\beta+\iota}|\geq |A_\beta|+|A_\iota|-1\geq |A_\beta|.\] Regardless, $|B_{\beta+\iota}|\geq |A_\beta|\geq \frac{p+1}{2}$.
    
    Thus, by Lemma \ref{Case2Backbone} we have  \[|2\hat{\;}A|\geq 3p+|A_0|+|B_{\beta+\iota}|\geq 4p+2\geq 4p.\] 
\end{proof}

\begin{lem}\label{Case2:L2S0}
    In Case 2, If $\ell=2$ and $s=0$ then $|2\hat{\;}A|\geq 4p$.
\end{lem}
\begin{proof}
    We let $S=\{\beta,\gamma\}$ such that \[p\geq |A_0|\geq |A_\beta|\geq |A_\gamma|\geq 1.\] By Theorem \ref{cdt}, there are at least $5$ distinct elements in the set \[2\{0,\beta,\gamma\}=\{0,\beta,\gamma,\beta+\gamma,2\beta,2\gamma\}.\] Since $0,\beta,\gamma$ are distinct by Lemma \ref{Case2Backbone} we then have that \[|2\hat{\;}A|\geq 3p+|B_{2\beta}|+|B_{2\gamma}|+|B_{\beta+\gamma}|-\max\{|B_{2\beta}|,|B_{2\gamma}|,|B_{\gamma+\beta}|\}.\] 
Note now that because $|A_0|+|A_\beta|+|A_\gamma|=2p+1$ and $|A_0|\leq p$ we have that $|A_\beta|+|A_\gamma|\geq p+1$, and so we have that \[|B_{\beta+\gamma}|\geq |A_{\beta}+A_\gamma|\geq |A_\beta|+|A_\gamma|-1\geq p.\] This then gives us that \[|2\hat{\;}A|\geq 3p+|B_{2\beta}|+|B_{2\gamma}|\geq 3p+|2\hat{\;}A_\beta|+|2\hat{\;}A_\gamma|\]
\[\geq 3p+\minz\{p,2|A_\beta|-3\}+\minz\{p,2|A_\gamma|-3\}.\] Thus, if $|2\hat{\;}A|\leq 4p-1$ then we must have that \[|2\hat{\;}A|\geq 3p+2|A_\beta|+2|A_\gamma|-6\geq 5p-4> 4p. \]
\end{proof}
For the case of $\ell=1$ we let $\beta$ be the unique element of $S$ such that $|A_0|+|A_\beta|-1\geq p$.
\begin{lem}
\label{Case2:L1reduce}
   In Case 2, If $\ell =1$ and $|2\hat{\;}A|\leq 4p-1$ then either 
    \begin{enumerate}
        \item $s=1$ or
        \item $s=2$ and $|A_\beta|=|A_\alpha|$. 
    \end{enumerate}
\end{lem}
\begin{proof}Note that $s\neq 0$ as we must have that $2\leq m=\ell+s=s+1$.

    We now observe that \[|2\hat{\;}A|\geq \sum_{i\in S_0}|B_i|\geq 2p+\sum_{i\in S\setminus\{\beta\}}|B_i|\geq 2p+\sum_{i\in S\setminus\{\beta\}}|A_0+A_i|\geq 2p+\sum_{i\in S\setminus\{\beta\}}(|A_0|+|A_i|-1).\]
    \[=2p+s(|A_0|-1)+(|A|-|A_0|-|A_\beta|)=4p+(s-1)(|A_0|-1)-|A_\beta|.\] 

    If $s\geq 3$ we have that \[|2\hat{\;}A|\geq 4p+2(|A_0|-1)-|A_\beta|\geq 4p+ |A_0|-2\geq 4p.\]

    If $s=2$, let us define $\delta=|A_0|-|A_\beta|$ and see that \[|2\hat{\;}A|\geq 4p+(|A_0|-1)-|A_\beta|\geq 4p-1+\delta,\] and so if $\delta\neq 0$ (or equivalently $|A_0|=|A_\beta|$) then $|2\hat{\;}A|\geq 4p$.

    Our claim now follows.
\end{proof}

\begin{lem}\label{Case2:L1S2}
    In Case 2, if $\ell=1$ and $s=2$ then $|2\hat{\;}A|\geq 4p$.
\end{lem}
\begin{proof}
In this case, consider the typical coset partition $A=A_0\cup A_\beta\cup A_\gamma\cup A_\delta$ with \[|A_0|\geq |A_\beta|\geq p+1-|A_0|>|A_\gamma|, |A_\delta|.\footnote{Unlike the case of $\ell=2$ and $s=1$ we may have that $|A_\gamma|<|A_\delta|$.}\] 
By Lemma \ref{Case2:L1reduce} we also have that $|A_0|=|A_\beta|$ and so \[|B_{2\beta}|\geq |2\hat{\;}A_\beta|\geq \minz\{2|A_\beta|-3,p\}\geq \minz\left\{2\frac{p+3}{2},p\right\}\geq\minz\{p,p\}=p.\]

If $2\beta\not\in \{\gamma,\delta\}$ then we have that $\{0,\beta,\gamma,\delta,2\beta\}$ are distinct and so we have that 
\[|2\hat{\;}A|\geq |B_0|+|B_\beta|+|B_{2\beta}|+|B_{\gamma}|+|B_\delta|\geq 3p+2|A_0|\geq 4p+3\geq 4p.\]

Assume WLOG then that $\delta=2\beta$. This then implies $\{0,\beta,2\beta,\gamma\}$ are pairwise  distinct. If $\delta+\beta$ is also pairwise distinct from these four then we similarly obtain\[|2\hat{\;}A|\geq |B_0|+|B_\beta|+|B_{2\beta}|+|B_{\gamma}|+|B_{\beta+\delta}|\geq 3p+2|A_0|\geq 4p+3\geq 4p.\]

Thus, if $|2\hat{\;}A|\leq 4p-1$ then $\delta+\beta\in\{0,\beta,2\beta,\gamma\}$, but clearly we cannot have $\delta+\beta=\beta$ or $\delta+\beta=2\beta$ We additionally see  that $\delta+\beta\neq 0$, as this would mean that $3\beta=0$ which cannot be as $p\geq 5$ and $\beta\neq 0$.

Thus, if $|2\hat{\;}A|\leq 4p-1$ then $\delta+\beta=\gamma$ which implies that $3\beta=\gamma$ and so \[S_0=\{0,\beta,2\beta,3\beta\}\] in which case (implied by the fact that $p\geq 5$) gives us
\[|2\hat{\;}A|\geq |B_0|+|B_{\beta}|+|B_{2\beta}|+|B_{3\beta}|+|B_{4\beta}|=3p+|B_{\delta}|+|B_{\beta+\delta}|\] \[\geq 3p +2|A_0|\geq 4p+2>4p.\]

% By Theorem \ref{cdt} there are at least $7$ distinct elements in the set $2S_0$. Via (\ref{CosetLBMain}) we may deduce that for $i\in 2S_0$ that
% \begin{equation}
%     |B_i|\geq\begin{cases}
%         p & i=0,\beta,2\beta\\
%         2|A_\gamma|-3 & i=2\gamma\\
%         2|A_\delta|-3 & i=2\delta\\
%         |A_\delta|+|A_\gamma|-1 & i=\delta+\gamma\\
%         |A_0| & \text{otherwise}\\
%     \end{cases}
% \end{equation}

% Since $\{0,\beta,\gamma,\delta\}$ are pairwise distinct, if for some $e\in\{\gamma,\delta\}$ is a member of the set $\{0,\beta,2\beta\}$ then we must have that $e=2\beta$, but by the fact that $\gamma\neq \delta$ we cannot have that both $\gamma=2\beta$ and $\delta=2\beta$, and thus for at least one $e_1\in \{\gamma,\delta\}$ the set $\{0,\beta,2\beta,e_1\}$ are four distinct elements. Letting $e_2\in\{\gamma,\delta\}$ but $e_2\neq e_1$ we see that the pair $\{e_2,e_2+\beta\}$ cannot be a subset of $\{0,\beta,2\beta,e_1\}.$ This is seen to be true as follows: If $\{e_2,e_2+\beta\}\subseteq \{0,\beta,2\beta,e_1\}$ then we must have that $e_2=2\beta$ and $e_2+\beta\in$
\end{proof}

\begin{lem}\label{Case2:L1S1}
    In Case 2, if $\ell=1$ then $s\neq 1$.
\end{lem}
\begin{proof}
    We have the coset partition $A=A_0\cup A_\beta\cup A_\gamma$ with \[|A_0|\geq |A_\beta|\geq p+1-|A_0|>|A_\gamma|\] and  \[|A_0|+|A_\beta|+|A_\gamma|=2p+1.\] Together, these imply that $|A_\beta|>p$, which cannot be.
\end{proof}
We now move to the final case: $\ell=0$.
\begin{lem}\label{Case2:L0}
  In Case 2, if $\ell=0$ then $|2\hat{\;}A|\geq 4p$.
\end{lem}
\begin{proof}
     If $\ell=0$ then it follows that for all non-zero $i$ we have $|A_i|+|A_0|-1\leq p-1$, and since $|A_0|\geq \frac{p+3}{2}$ we have that \[|A_i|\leq \frac{p-3}{2}.\] 

    Additionally, observe that for any $j\in S$ we have that \[2p+1=|A|=\sum_{i\in S_0}|A_i|=(|A_0|+|A_j|)+\sum_{i\in S\setminus \{j\}}|A_i|\leq p+(s-1)\frac{p-3}{2},\] and so we have that 
    \[s\geq 1+\frac{2p+2}{p-3}=3+\frac{8}{p-3}>3,\] and so we must have $s\geq 4$.

    We now use (\ref{basiclb}), Lemma \ref{Case2:2A0}, and (\ref{CosetLBMain}) to get \[|2\hat{\;}A|=\sum_{i\in T_0}|B_i|\geq |2\hat{\;}A_0|+ \sum_{i\in S}|A_0+A_i|\geq p+\sum_{i\in S}(|A_0|+|A_i|-1)=p+s(|A_0|-1)+(|A|-|A_0|)\]
    \[=3p+(s-1)(|A_0|-1).\]
    We now recall that $s\geq 4$ and $|A_0|\geq \frac{p+3}{2}$ and so we have that \[|2\hat{\;}A|\geq \frac{9p+9}{2}> 9p/2>4p.\]
\end{proof}

\begin{lem}
    In Case 2, we have that $|2\hat{\;}A|\geq 4p$.
\end{lem}
\begin{proof}
    If $\ell\geq 3$ then our claim follows from Lemma \ref{Case2:lg3}. If $\ell=2$ then our claim follows from Lemma \ref{Case2:L2S2}, Lemma \ref{Case2:L2S1}, and Lemma \ref{Case2:L2S0}. If $\ell=1$ then our claim follows from Lemma \ref{Case2:L1reduce}, Lemma \ref{Case2:L1S2}, and Lemma \ref{Case2:L1S1}. Lastly, if $\ell=0$ then our claim follows from Lemma \ref{Case2:L0}.
\end{proof}

\section{Conclusion}

With this, we have that regardless of $m$ or $|A_0|$: $|2\hat{\;}A|\geq 4p$ for all $A\subseteq \Z_p^2$ for $p\geq 11$ and $|A|=2p+1$ and so Theorem \ref{maintheorem} is proven for all $p$ except $p=5,7$. Using \cite{SumsetGenerator} with a powerful enough computer verifies the theorem for these two values of $p$, completing the proof of Theorem \ref{maintheorem}.

While this result is a major step forward, the author advises caution for a reader who wishes to generalize this result using the methods in this paper. There are two potential directions for generalization. The first is relaxing the condition of $G\cong \Z_p^2$ to $G\cong \Z_p^r$ for some $r$. While Case 1B (and to a lesser extent Case 2) seem capable make this generalization with only a few minor issues, Case 1A's reduction to the other two cases relies explicitly on both $H\cong \Z_p$ and $G/H\cong \Z_p$ which is only possible in the case of $G\cong \Z_p^2$. In order to prove that $\rho(\Z_p^r,2p+1)=4p$ is true for a sufficiently large $p$, a new method must be developed for Case 1A.

Regardless, the author believes that Theorem \ref{maintheorem} generalizes in its entirety. Specifically:
\begin{conj}
    If $p\geq 5$ is prime then \[\rho(\Z_p^r,2p+1)=4p.\]
\end{conj}

The second way that the results of this paper can be generalized is by determining $\rho(\Z_p^2,kp+1)$ for $k\geq 3$. Like before, Case 1B seems to generalize rather nicely, and Case 1A also does not appear to have any outstanding issues regarding its generalization (except perhaps, a stricter lower bound on when Lemma \ref{Case1A:Reduction} reduces the problem to Cases 1B and 2). The problem occurs when examining Case 2. Here, to prove that Case 2 cannot provide a counterexample to $\rho(\Z_p^2,2p+1)=4p$ we considered each value of $\ell$, one at a time. However, if one were to go out and prove, say $\rho(\Z_p^2,3p+1)=6p$ they would need to consider $\ell\leq 4$ if they wanted to directly adapt the methods used in this paper. And if one wishes to consider the general case of proving that $\rho(\Z_p^2,kp+1)=2kp$ then they will need to consider every case when $\ell\leq 2kp-2$ which will require a less ``brute force'' approach than what is used in Section \ref{Section2}.

Another natural question following the results of this paper is to solve the corresponding ``inverse problem'' of Theorem \ref{maintheorem}, ie. the problem of classifying all $2p+1$-sets $A\subseteq\Z_p^2$ such that $|2\hat{\;}A|=\rho(\Z_p^2,2p+1)=4p$. The equivalent problem for sets of size $p+1$ was solved in \cite{EliKer} as follows.
\begin{thm}[Eliahou and Kervaire \cite{EliKer}]
    For prime $p\geq 5$ and $A\subseteq \Z_p^r$ if $|A|=p+1$ and $|2\hat{\;}A|=2p$ then there exists an order $p$ subgroup $Z< \Z_p^r$ such that $A$ is the union of a coset of $Z$ and a single element outside of said coset.
\end{thm}
We believe that our case is rather similar and offer the following conjecture:
\begin{conj}
    For prime $p\geq 7$ then for $A\subseteq \Z_p^2$ if $|A|=2p+1$ and $|2\hat{\;}A|=4p$ then there exists an order $p$ subgroup $Z<\Z_p^2$ with canonical homomorphism $\phi: \Z_p^2\to \Z_p^2/Z$ such that $\phi(A)$ is an arithmetic progression of length three and that there is a unique element $a\in A$ such that $A\setminus \{a\}$ is the union of two cosets of $Z$.  
\end{conj}

% \begin{itemize}
%     \item Perhaps talk about why this proof doesn't work for $\rho(\Z_p^r,2p+1)$?
%     \begin{itemize}
%         \item Case 1B doesn't adapt well.
%     \end{itemize}
%         \item Perhaps talk about why this proof doesn't work for $\rho(\Z_p^r,3p+1)$? or $\rho(\Z_p^2,kp+1)$ for $k>3$?
%     \begin{itemize}
%         \item Case 2 doesn't adapt well.
%     \end{itemize}
% \end{itemize}

\bibliographystyle{plain} % We choose the "plain" reference style
\bibliography{rho} % Entries are in the refs.bib file

@book{Baj:2018a,
      title={Additive Combinatorics: A Menu of Research Problems}, 
      author={Bela Bajnok},
      publisher={CRC Press, Boca Raton, FL}, 
    doi={https://doi.org/10.1201/9781351137621},
    series={Discrete Mathematics and Its
Applications.},
    year={2018},
}

@article{EliKer,
title = {Restricted sums of sets of cardinality $1+p$ in a vector space over $\mathbb{F}_p$},
journal = {Discrete Mathematics},
volume = {235},
number = {1},
pages = {199-213},
year = {2001},
note = {Chech and Slovak 3},
issn = {0012-365X},
doi = {https://doi.org/10.1016/S0012-365X(00)00273-9},
url = {https://www.sciencedirect.com/science/article/pii/S0012365X00002739},
author = {Shalom Eliahou and Michel Kervaire}
}

@article{DDSandH,
author = {Da Silva, J. A. Dias and Hamidoune, Y. O.},
title = {Cyclic Spaces for Grassmann Derivatives and Additive Theory},
journal = {Bulletin of the London Mathematical Society},
volume = {26},
number = {2},
pages = {140-146},
doi = {https://doi.org/10.1112/blms/26.2.140},
url = {https://londmathsoc.onlinelibrary.wiley.com/doi/abs/10.1112/blms/26.2.140},
eprint = {https://londmathsoc.onlinelibrary.wiley.com/doi/pdf/10.1112/blms/26.2.140},
year = {1994}
}

@article{Vosper,
author = {Vosper, A. G.},
title = {The Critical Pairs of Subsets of a Group of Prime Order},
journal = {Journal of the London Mathematical Society},
volume = {s1-31},
number = {2},
pages = {200-205},
doi = {https://doi.org/10.1112/jlms/s1-31.2.200},
url = {https://londmathsoc.onlinelibrary.wiley.com/doi/abs/10.1112/jlms/s1-31.2.200},
eprint = {https://londmathsoc.onlinelibrary.wiley.com/doi/pdf/10.1112/jlms/s1-31.2.200},
year = {1956}
}

@misc{SumsetGenerator,
author = {Llano, John},
title = {Sumset Generator},
 howpublished={\url{http://cs.gettysburg.edu/~llanjo01/addcomb/sumsets.html}},
year = {2024}
}

@article{SmallSums,
title = {Optimally small sumsets in finite abelian groups},
journal = {Journal of Number Theory},
volume = {101},
number = {2},
pages = {338-348},
year = {2003},
issn = {0022-314X},
doi = {https://doi.org/10.1016/S0022-314X(03)00060-X},
url = {https://www.sciencedirect.com/science/article/pii/S0022314X0300060X},
author = {Shalom Eliahou and Michel Kervaire and Alain Plagne}
}

@article{EKfirst,
title = {Sumsets in Vector Spaces over Finite Fields},
journal = {Journal of Number Theory},
volume = {71},
number = {1},
pages = {12-39},
year = {1998},
issn = {0022-314X},
doi = {https://doi.org/10.1006/jnth.1998.2235},
url = {https://www.sciencedirect.com/science/article/pii/S0022314X98922351},
author = {Shalom Eliahou and Michel Kervaire}
}
\end{document}